# A Fundamental Theorem of Powerful Set-Valued for F-Rough Ring

Faraj.A.Abdunabi[1]
[1]Mathematics department
University of Ajdabyia
Benghazi, Libya

Ahmed shletiet[2]
[2]Mathematics department
University of Ajdabyia
Benghazi, Libya

**Abstract:-** In this paper, we introduce the upper and lower approximations on the invers set-valued mapping and the approximations an established on a powerful set valued homomorphism from a ring $R_1$ to power sets of a ring $R_2$. Moreover, the properties of lower and upper approximations of a powerful set valued are studies. In addition, we will give a proof of the theorem of isomorphism over approximations F-rough ring as new result. However, we will prove the kernel of the powerful set-valued homomorphism is a subring of $R_1$. Our result is introduce the first isomorphism theorem of ring as generalized the concept of the set valued mappings.

## I. INTRODUCTION

The rough set theory has been introduced by Pawlak 1982[1], as the new tool to incomplete information system. Many researchers develop this theory in many areas. Substantively, the rough set an established on two concepts of approximations (lower and upper). [2],R. Biswas and S. Nanda Show  rough groups and subgroups. N. Kuroki [3] consider the rough the notation of ideal in a semi-grousp. B. Davvaz have studied roughness in ring [4]. The concepts of rough prime ideals and rough primary ideals in a ring has introduced by O. Kazanci and B. Davvaz [5],. V. Selvan and G. Senthil Kumar [6] consider notation of the roughness ideals on a semi-ring. In [7], B. Davvaz given the set valued homomorphism and study the T-rough sets in a group. In [8], B. Davvaz and others has generalized the concepts of upper and lower approximations established on a ring by the set valued homomorphism of rings. The properties of T- rough sets in commutative rings has studied by.S. B. Hosseini, N. Jafarzadeh, and A. Gholami[9]. However, other researchers have been interested of the Set-valued maps [10],[11]. Set-valued maps have used in many areas such as Economics [12]. In [13], G. SenthilKumarthis gives the proof the fundamental set-valued homomorphism group theorem. Our work, we study the powerful set-valued homomorphism established on ring and some their propertiess. Moreover, we will show that the kernel of the powerful set-valued homomorphism is a subring of $R_1$. Our result is introducethe fundamental isomorphismtheorem of ring as generalized the concept of the set valued mappings.

## II. PRELIMINARIES

We will recall the concept of the approximation (lower, lower) based on a set valued mapping for more information and proofs we can see [7]. Ina addition, we introduce the set valued homomorphism.

***Definition 2-1***: Suppose $U \neq \varnothing$. Let ~ be an equivalence relation on $U$. Let $R: 2^U \to 2^U \times 2^U$ where $2^U$ is the set of all non-empty subsets of $U$. A pair $(U, \sim)$ is called an approximation space and the upper rough approximation of $X$ is $\overline{R(X)} = \{x \in U: [x]_\sim \subseteq X\}$ and $\underline{R(X)} = \{x \in U: [x]_\sim \cap X \neq \varnothing\}$ is th lower rough approximation of $X$.

***Definition 2-2***: Suppose that $X, Y$ any two nonempty sets. A *set-valued map* or *multivalued map* $F: X \to 2^Y$ from $X$ to $Y$ is a map that related to any $x \in X$ a subset of $F$ of $Y$. The set called the image of $x$ under $F$ by $F(x)$ and we define the domain of $F$ by $D_F = \{x \in X: F(x) \neq \varnothing\}$. The image of $F$ is a subset of $Y$ defined by $\text{Im}(F) = \bigcup_{x \in X} F(x) = \bigcup_{x \in D_F} F$.

***Definition 2.3:*** Suppose that $X$; $Y$ any non-empty sets. Let $M \subseteq Y$ and $F: X \to 2^Y$ be a set-valued mapping. We called the $\underline{F(M)} = \{x \in X: F(x) \subseteq M\}$; $\overline{F(M)} = \{x \in X: F(x) \cap M \neq \varnothing\}$ respectively and $(\underline{F(M)}, \overline{F(M)})$ is called $F$- rough set of $F$.

***Example 2-1***: let $X = \{1, 2, 3, 4, 5, 6\}$ and let $F: X \to 2^X$ where $\forall x \in X$, $F(1) = \{1\}$, $F(2) = \{1, 3\}$, $F(3) = \{3,4\}$, $F(4) = \{4\}$, $F(5) = \{1,6\}$, $F(6) = \{1, 5, 6\}$. Let $A = \{1, 3, 5\}$, then$\underline{F(M)} = \{1, 2\}$, and $\overline{F(A)} = \{1, 2, 3, 5, 6\}$, B(A)$\neq \varnothing$, is rough. $\text{Im}(F) = \bigcup_{x \in X} F(x) = \{1,3,4,5\}$ Now, Let $B = \{2,4,6\}$= then$\underline{F(B)} = \{4\}$, and $\overline{F(B)} = \{3,4, 5, 6\}$, $B(B) \neq \varnothing$, is rough.

***Preposition 2-1:*** Suppose that $X$, $Y$ are non-empty sets. Let $A, B \subseteq Y$. If $F: X \to 2^Y$ be a set-valued mapping , then:
1)- $\overline{F(A) \cap F(B)} = \overline{F(A)} \cap \overline{F(B)}$;
2) $\overline{F(\varnothing)} = \varnothing$;
3)- $\overline{F(A) \cup F(B)} = \overline{F(A)} \cup \overline{F(B)}$;

4)- F(A) ∪ F(B)= $\overline{F(A)} \cup \overline{F(B)}$;
5)- $\overline{F(X)}$ = $\underline{F(A)} \cup \underline{F(B)}$

*Proof:* By using the definitions upper and lower approximations

*Example 2-2*: let $X = \{a, b, c, d, e, f\}$ and let $F : X \to 2^X$. If $F_1(a) = \{a\}$, $F_1(b) = \{a, c\}$, $F_1(c)= \{c,d\}$, $F_1(d) = \{a,d\}$, $F_1(e) = \{a,f\}$, $F_1(f) = \{a, e, f\}$. And $F_2(a) = \{a\}$, $F_2(b) = \{a, b\}$, $F_2(c) = \{c\}$, $F_2(d) = \{d\}$, $F_2(e) = F_2(f) =\{a,e,f\}$ and $A = \{a, c, e\}$, then $\underline{F_1(A)} = \{a, b\}$, and $\overline{F_1(A)} = \{a, b, c, e, f\}$, $\underline{F_2(A)} = \{a, c\}$, $\overline{F_2(A)} = \{a, b, c, e, d, f\}$,

*Theorem 2-1* Let $F_1; F_2: X \to 2^X$ be set-valued map such that $(F_1 \cap F_2)x \neq \emptyset$ ; $\forall x \in X$ and $A \subseteq X$. Then
1- $\overline{F_1(A) \cup F_2(A)} = \overline{F_1(A)} \cup \overline{F_2(A)}$;
2- $\overline{F_1(A) \cap F_2(A)} \subseteq \overline{F_1(A)} \cap \overline{F_2(A)}$;
3- $\underline{F_1(A) \cap F_2(A)} = \underline{F_1(A)} \cap \underline{F_2(A)}$;
4- $\underline{F_1(A) \cup F_2(A)} = \underline{F_1(A)} \cup \underline{F_2(A)}$.

*Proof*: By using, the definitions 2-3.

*Definition 2.4*. Let $2^X$ be the set of all subsets of a non-empty $X$. If $S_1, S_2 \in 2^X$, then we define $S_1+S_2 =\{x \in X;$ either $X \in S_1$ or $x \in S_2$, but not in both$\}$ and $S_1 * S_2 = S_1 \cap S_2$ are called sum and product of $S_1$ and $S_2$ respectively.

*Preposition 2-2* Let $2^X$ be the power of all subsets of a non-empty $X$ with sum and product of $S_1$ and $S_2$. Then $(2^X, *, +)$ is a commutative ring.

Note that, the empty set $\emptyset$ is the identity of $+$ and the set $X$ is the identity of $*$. Therefore, we called $2^X$ is Ring of Subsets of $X$.

## III.　　MAIN RESULTS

We introduce the concepts of the invers set-valued map. Suppose $X$ and $Y$ are two nonempty sets.

*Definition 3-1*: Suppose that $F$ is a set-valued map from $X$ to $Y$, we call $F^{-1}$ the inverse of $F$ and we write as: $F^{-1}(y)= \{ x \in X: y \in F(x)\}$, $\forall y \in Y$. If $B \subseteq Y$, then the upper inverse image is $\overline{F^{-1}}$ $=\{x \in X : F(x) \cap B \neq \emptyset\}$ and lower inverse image is $\underline{F^{-1}}$ $=\{x \in X: F(x) \subseteq B\}$ and the boundary is B($F^{-1}$)= $\overline{F^{-1}}$ - $\underline{F^{-1}}$. If B($F^{-1}$)$\neq \emptyset$ is called roughly invers a set-valued map.

*Example 3-1*: Consider the example 1-1, let $X = \{a, b, c, d, e, f\}$. Suppose $F : X \to 2^X$ is set value map where $\forall x \in X$, $F(a)= \{a\}$, $F(b)=\{a, c\}$, $F(c) = \{c,d\}$, $F(d)=\{d\}$, $F(e)=\{a,f\}$, $F(f) = \{a, e, f\}$. Let $B_1 = \{a, d, e\}$, then $\overline{F^{-1}(B_1)} = \{a, d\}$, and $\underline{F^{-1}(B_1)} = \{a, b, c, d, e, f\}$, $B(B_1) \neq \emptyset$, is rough. Let $B_2=\{b,e,f\}=$ then $\overline{F^{-1}(B_2)} = \{d\}$, and $\underline{F^{-1}(B_2)} = \{ c, d, e, f\}$, $B(B_2) \neq \emptyset$, is rough.

*Proposition 3-1*: Let $X, Y$ be non-empty sets and $B_1, B_2 \subseteq Y$. If $F:X \to 2^Y$ be a set-valued mapping, then:
1)- $\overline{F^{-1}(B_1) \cap F^{-1}(B_2)} = \overline{F^{-1}(B_1)} \cap \overline{F^{-1}(B_2)}$;
2)- $\overline{F^{-1}(B_1) \cup F^{-1}(B_2)} = \overline{F^{-1}(B_1)} \cup \overline{F^{-1}(B_2)}$;
3)- $F^{-1}(B_1) \cup F^{-1}(B_2) = \underline{F^{-1}(B_1)} \cup \underline{F^{-1}(B_2)}$

*Proof*: it is explicit.

*Proposition 3-2*: Let $F_1; F_2 : X \to 2^Y$ be set-valued map such that $(F_1 \cap F_2)(x) \neq \emptyset$ ; $\forall x \in X$ and $B \subseteq Y$. Then
1)- $\overline{(F_1 \cup F_2)^-}(B) = \overline{F_1}^-(B) \cup \overline{F_2}^{-1}(B)$;
2)- $\overline{(F_1 \cap F_2)^{-1}}(B) \subseteq \overline{(F_1 \cap F_2)^{-1}}(B)$;
3)- $\underline{(F_1(B) \cup F_2(B))^{-1}} = \underline{(F_1(B))^{-1}} \cap \underline{(F_2(B))^{-1}}$ ;
4)- $\underline{(F_1(B) \cap F(B))^{-1}} = \underline{(F_1(B)^{-1}} \cup (F_2(B))^{-1}$.

*Proof*: it is explicit

*Example 3-2:* Suppose that $U = \{1, 2, 3, 4, 5, 6\}$. let $F : U \to 2^U$ where $\forall x \in X$, $F(1) = \{1\}$, $F(2) = \{1, 3\}$, $F(3) = \{3,4\}$, $F(4) = \{4\}$, $F(5) = \{1,6\}$, $F(6) = \{1, 5, 6\}$. (1) If $A = \{1, 3, 5\}$, then $\overline{F(A)} = \{1, 2, 3\}$, and $\underline{F(A)} = \{1, 2, 3, 4, 5, 6\}$, $B(A) \neq \emptyset$, is rough.

## IV.　　ISOMORPHISM THEOREM FOR F-ROUGH RINGS.

*Definition 4.1:* suppose that $R$ is a ring. Let $\sim$ be a conformity of $R$, that is, $\sim$ is an equivalence relation on $R$ such that $(x, y) \in \sim$ implies $(xa, ya) \in \sim$ and $(ax, ay) \in \sim$ for all $a \in R$. We denote by $[x]_\sim$ the $\sim$ conformity class containing the element $x \in R$.

*Remark 4-1:* Let $\sim$ be a conformity on a ring $R$. Define $F:R \to 2^R$ by $F(x)=[x]_\sim$, $\forall x \in R$. Then $F$ is a set-valued homomorphism.

*Definition 4.2*: Suppose that $R_1$ and $R_2$ are two rings. Let $F:R_1 \to 2^{R_2}$ be a set-valued homomorphism and $S$ be subring of $R_2$. If $\underline{F}$ is subring of $R_1$, then $S$ is called the lower F-rough subring of $R_2$ and if $\overline{F}$, is subring of $R_1$ then S is called the upper F-rough subring of $R_2$. If $\overline{F}$, and $\underline{F}$ are subring of $R_1$, then we call ($\overline{F}$, $\underline{F}$) a F-rough subring.

*Definition 4.3*: Suppose that $R_1$, $R_2$ are rings. *A powerful set–valued set-valued homomorphism* is a mapping from $R_1$ into $2^{R_2}$ that preserves the ring operations $(*, +)$, that is, $\forall$ l $x$ and $y \in R_1$ such that:
1) $F(x*y)= \{ab: a \in F(x_1), b \in F(x_2)\}$; and
2) $F(x+y) = F(x)+F(y)$;
3) $(F(x))^{-1} = \{ a^{-1}: a \in F(x)\} = F(a^{-1})$;
4)- $F(x*y)=F(x)F(y)$.

***Remark 4-2***: If the powerful set–valued set-valued homomorphism is one to one, then we called an *epimorphism*. If it is onto and one-to-one, then we called an *isomorphism*.

***Proposition 4.2:*** Let ~ be a conformity on a ring $R$. If $x,y \in R$ then:
1)- $[x*y]_\sim = [x]_\sim * [y]_\sim$ ;
2)- $[x+y]_\sim = [x]_\sim + [y]_\sim$ ;
3)- $([x]_\sim)^{-1} = \{a^{-1}: a \in [x]_\sim)\} = \{[x]_\sim\}^{-1}$.
*Proof*: from definition 4-1

***Proposition 4.2***: Suppose that $R_1$ and $R_2$ are two rings. Let $F: R_1 \to R_2$ powerful set–valued set-valued homomorphism from then:
1. $F(e_{R1}) = e_{R2}$;
2. For any $r \in R_1$, $F(-r) = -F(r)$. If $r$ is a unit, then $F(r^{-1}) = F(r)^{-1}$.
3. $F(R_1)$ is a subring of $R_2$.

***Theorem 4.1*** Suppose that $R_1$, $R_2$ are rings. Suppose $F(e_{R1}) \in S$ where S is a subring of $R_2$. If $F(e_{R1})$ and $F: R_1 \to 2^{R_2}$ be a powerful set-valued homomorphism, then $\underline{F}$ is a subring of $R_1$.

*Proof*: Suppose that $S$ is a subring of $R_2$. We need verify that the conditions of subring are satisfied by $\underline{F}$. We have $(S,+)$ is commutative subgroup and $(S, *)$ is semi-group. Suppose that $S$ be a subring of $R_2$ containing $F(e_{R1})$. Suppose that $e_{R2}$ be the identity of the ring $R_2$. The $\ker(F) = F^{-1}\{e_{R1}\} = \{r \in R_1 ; F(r) = e_{R2}\}$ is the preimage of $e_{R2}$. We have $\underline{F}^{-1} = \{e_{R1} \in R_1 : F(e_{R1}) \subseteq S\}$, then $e_{R1} \in \underline{F}$. So, $\underline{F}$ is nonempty. Let $x_1$ and $y_2$ be arbitrary elements of $\underline{F}$. Then there exist elements $x, y$ such that $F(x) = x_1$ and $F(y) = y_1$. Since $S$ is a subring of $R_2$, then $x+y$ and $x*y$ are in $S$. Therefore, $F(x+y) = F(x)+F(y) = x_1+y_2$ and $F(x*y) = F(x)*F(y) = x_1*y_1$ are $\underline{F}$. Now, Since $-x$ in $S$, $F(-x) = -F(x)$ (from Proposition 4.2). So, $\underline{F}$ is closed under + and *. Then $\underline{F}$ is a subring of $R_1$.

***Theorem 4.2***: Let $S$ and $R$ be two rings, Let $\rho: S \to R$ be subjective homomorphism and $F_2: S \to 2^R$ be a set-valued homomorphism. If $F_1(x) = \{s_1 \in S: \rho(s_1) \in F_2(\rho(x))\}$, $\forall x \in S$ and $A \subseteq R$, then $\rho(\overline{F_1}) = \overline{F_2(\rho}$.

*Proof*: suppose $y \in \rho(\overline{F_1})$. We have $x \in \overline{F_1}$ such that $y = \rho(x)$. Since $x \in \overline{F_1}$, then we get the exits element in $F_1(x)$ and $A$ that mean $F_1(x) \cap A \neq \emptyset$. Hence, there exists $a \in A$ with $a \in F_1(x)$. Thus $\rho(a) \in F_2(\rho(x))$ and $\rho(a) \in \rho(A) \Rightarrow F_2(\rho(x)) \cap \rho(A) \neq \emptyset$. Therefore, $y = \rho(x) \in \overline{F_2(\rho}$. Conversely, if $y \in \overline{F_2(((}$, then since $\rho$ is onto. We give us there exists $x \in R_1$ such that $y = \rho(x) \in \overline{F_2(\rho}$. That is, $F_2(\rho(x)) \cap \rho(A) \neq \emptyset$. So, there exists $z \in F_2(\rho(x)) \cap \rho(A)$. Thus, $z = \rho(a)$ for some $a \in A$. So, $z = \rho(a) \in F_2(\rho(x))$. If $a \in F_1(x) \cap A$, then $F_1(x) \cap A \neq \emptyset$. Therefore, we get $x \in \overline{F_1(\rho}$. Also,

$y \in \rho(x) \in \rho(\overline{F_1})$. So, $\overline{F_2(\rho} \subseteq$ then $\rho(\overline{F_1})$. We conclude that $\rho(\overline{F_1}) = \overline{F_2(f(}$ as we required.

***Theorem 4.3:*** Let $R_1$, $R_2$ be rings. Suppose $\rho: R_1 \to R_2$ be an epimorphism $F_2: R_1 \to 2^{R_2}$ be a powerful set-valued homomorphism. If $\rho$ is injective mapping and $F_1(x) = \{r_1 \in R_1: \rho(r_1) \in F_2(\rho(x))\}$ for all $x \in R_1$, then $F_1: R_1 \to 2^{R_1}$ is a powerful set valued homomorphism.

*Proof*: Suppose $s \in F_1(xy)$. We have $\rho(s) \in F_2(\rho(xy)) = F_2(\rho(x)\rho(y)) = F_2(\rho(x))*F_2(\rho(y))$. So, we get $\rho(s) = ab$, for some $a \in F_2(\rho(x))$, $b \in F_2(\rho(y))$. Because, $\rho$ is subjective then there exists $c,d \in R_1$ with $\rho(c) = a$, $\rho(d) = b$. We get, $\rho(s) = \rho(c)\rho(d)$, $c \in F_1(x)$ and $d \in F_1(y)$. Therefore, $s = cd \Rightarrow s \in F_1(x)F_2(y)$. We get, $F_1(xy) \subseteq F_1(x)F_2(y)$.

On the other side, we suppose that $z \in F_1(x)F_1(y)$. Then we get $z = cd$ for some $c \in F_1(x)$ and $d \in F_1(y)$ and therefore, $(c) \in F_2(\rho(x), \rho(d) \in F_2(\rho(y))$. So, we get $\rho(cd) = \rho(c)\rho(d) \in F_2(\rho(x))F_2(\rho(y)) = F_2(\rho(x)(\rho(y)) = F_2(\rho(xy)) \Rightarrow z = cd \in F_1(xy)$. So, $F_1(x)F_2(y) \subseteq F_1(xy)$. So we get, $F_1(xy) = F_1(x)F_1(y)$.

Now, we assume that $s \in F_1(x+y)$. Then $\rho(u) \in F_2(\rho(x+y)) = F_2(\rho(x) + \rho(y)) = F_2(\rho(x)) + F_2(\rho(y))$. We get $\rho(s) = a+b$, for some $a \in F_2(\rho(x))$ and $b \in F_2(\rho(y))$. We know that $\rho$ is surjective, then we have $c, d \in R_1$ such that $\rho(c) = a$ and $\rho(d) = b$. subsequently, we have $\rho(s) = \rho(c) + \rho(d)$, $c \in F_1(x)$, and $d \in F_1(y)$. Therefore, $s = c + d \Rightarrow u \in F_1(x) + F_2(y)$. So, $F_1(x+y) \subseteq F_1(x) + F_2(y)$.

On the other side, assume that $z \in F_1(x) + F_1(y)$. Then $z = c+d$ for some $c \in F_1(x)$ and $d \in F_1(y)$ and so, $\rho(c) \in F_2(\rho(x)), \rho(d) \in F_2(\rho(y))$. So, $\rho(c+d) = \rho(c) + \rho(d) \in F_2(\rho(x)) + F_2(\rho(y)) = F_2(\rho(x) + (\rho(y)) = F_2(\rho(x+y)) \Rightarrow$, $z = c + d \in F_1(xy)$. We get $F_1(x) + F_2(y) \subseteq F_1(xy)$. So, we conclude that $F_1(x+y) = F_1(x) + F_1(y)$.

Finally, we have $x \in F_1(x^{-1}) \Leftrightarrow \rho(c) \in F_2(\rho(x^{-1})) \Leftrightarrow f(c) \in F^2(\rho(x))^{-1}) \Leftrightarrow \rho(c) \in F_2(\rho(x)))^{-1} \Leftrightarrow (\rho(c))e \in F_2(\rho(x)) \Leftrightarrow \rho(c^{-1}) \in F_2(\rho(x)) \Leftrightarrow c^{-1} \in F_1(x) \Leftrightarrow c \in (F_1(x))^{-1}$. Therefore, $F_1(x^{-1}) = (F_1(x))^{-1}$, for all $x \in R_1$. We get $F_1$ is a powerful set-valued homomorphism as we required.

***Definition 4-4***: Suppose that $R_1, R_2$ are rings. Let $F: R_1 \to 2^{R_2}$ be a powerful set-valued homomorphism. a kernel of $F$ is $Ker(F) = \{x \in R_1: F(x) = F(e_1)\}$, where $e_1$ is the identity element of $R_1$.

***Theorem 4.4*** Suppose that $R_1$ and $R_2$ are two rings. If $F: R_1 \to 2^{R_2}$ be a powerful set-valued homomorphism, then $Ker(F)$ is a subring of $R_1$.

*Proof*: we assume that $x, y \in Ker(F)$. From definition of $ker(F)$, $F(x) = F(e_1)$ and $F(y) = F(e_1)$. We have $F(xy) = F(x)F(y) =$

$F(e_1)F(e_1)=F(e_1)$ (*Since F is a powerful set valued homomorphism*).

Therefore, $xy \in Ker(F)$ and $F(x^{-1})=(F(x))^{-1}= (F(e_1))^{-1}=F(e_1^{-1})=F(e_1) \Rightarrow x^{-1} \in Ker(F)$.

Now, suppose $x+y \in Ker(F)$. Then $F(x+y) =F(x)+F(y) =F(e_1)+F(e_1) \in Ker(F)$. Therefore, $x+y \in Ker(F)$. we get, $Ker(F)$ is a subring of $R_1$ as we required.

## V. THE FUNDAMENTAL THEOREM OF RING HOMOMORPHISM.

Here if we have a ring $R$ is a ring and the ideal $I$ and $S$ subring of $R$, then the quotient of $R$ by $I$ is the set $R/I$ of equivalence classes $a+I=\{a+i:i \in I\}$ with the two operations$(+,*)$. We define $(a+I), (b+I) \in R/I$ by: $(a+I)+(b+I)=(a+b)+I$; and $(a+I)*(b+I)=(a*b)+I$ $\forall$ $(a+I), (b+I) \in R/I$. If $F: R \to 2^{S2}$ is a powerful set-valued homomorphism from a ring to $2^S$, then $ker(F)$ is a subring of $R$ *(Theorem 4.4)*. We would ultimately like $ker(F)$ to be an ideal of $R$ to define the quotient $R/ker(F)$.

Now, suppose that $e \in R$ with $a*e=a$ and $e*a=a$ $\forall$ $a \in R$, we need to find identity of operation $+$. By definition of $ker(F)$ in general cannot be a subring of $R$ because if e is the identity of $R$, then by definition of a ring homomorphism, $F(e)$ is mapped to the multiplicative identity of $S$ and not to the additive identity of $S$. To establish a fundamental theorem of ring homomorphisms, we make a small exception in not requiring that $ker(F)$ is an ideal for the quotient $R/ker(F)$ to be defined.

**Theorem5-1:**(*The Fundamental Theorem of Ring Homomorphisms*): Let $R$ and $S$ be two rings with ring. Let $F: R_1 \to 2^{R2}$ be a powerful set-valued homomorphism.. Then $R/ker(F) \cong F(R)$.
*Proof*: Let $ker(F)$ be the kernel of $F$. let $\phi:R/Ker(F) \to F(R)$ and $\forall$ $(a+Ker(F)) \in R$ is $\phi(a+Ker(F))=F(a)$. We show that $\phi$ is well-defined. For $a+Ker(F)=$
$b+Ker(F) \Rightarrow a=b+ker(F)$ for some $k \in Ker(F)$. So, $\phi(a+Ker(F))=\phi((b+k)+Ker(F))=\phi$
$((b+Ker(F))+Ker(F))=\phi(b+Ker(F))$. Now, If $(a+Ker(F))$ and $(b+Ker(F)) \in R/Ker(F)$, then $\phi((a+Ker(F))+(b+Ker(F)))= \phi((a+b)+Ker(F))$
$=F(a+b)=F(a)+F(b)=\phi(a+Ker(F))+\phi(b+Ker(F));$
$\phi((a+Ker(F))*(b+Ker(F)))=\phi((a*b)+Ker(F))=F(a*b)=F(a)*\phi(b)= \phi (a+Ker(F))*\phi (b+Ker(F))$.
Finally, we need to show $\phi$ is bijective. Let $(a+Ker(F)),(b+Ker(F)) \in R/Ker(F)$ and suppose that $\phi (a+Ker(F))= \phi (b+Ker(F))$. Then $F(a)=F(b)$. So $F(a-b)=e$. So $a-b \in Ker(F)$. So $a+Ker(F)=b+Ker(F)$. Hence $\phi$ is injective. Now, $\forall a \in F(R)$ we have that $(a+I) \in R/Ker(F)$ is such that $\phi (a+Ker(F))=a$. So $\phi$ is surjective. That mean $\phi$ is bijective. So, $\phi$ is an isomorphism from $R/Ker(F)$ to $F(R)$, that is: $R/ker(F) \cong F(R)$ as required.

## VI. CONCLUSION

Theoretically, rough set based on the upper and lower approximations as an equivalence relation. In this paper, we introduce the upper and lower approximations on the invers set-valued mapping and the approximations an established on a powerful set valued homomorphism from a ring $R_1$ to the power sets of a ring $R_2$. Moreover, the properties of lower and upper approximations of a powerful set valued are studies. In addition, we will give a proof of the isomorphism theorem for lower and upper F-rough ring as new result. However, we will prove the kernel of the powerful set-valued homomorphism is a subring of $R_1$. Our result is introduce the first isomorphism theorem of ring as generalized the concept of the set valued mappings. We hope that it will be useful in some applications in the future.

## ACKNOWLEDGEMENT

The author would be very happy and would welcome any suggestions for improving this paper.

...